\documentclass{amsart}

\usepackage{amsmath,amsthm,amsfonts,amscd,amssymb}
\usepackage[pdfborder={0 0 0}]{hyperref}
\usepackage{enumerate}

\newcommand{\ZZ}{\mathbb{Z}}

\newtheorem{thm}{Theorem}[section]

\newtheorem{lem}[thm]{Lemma}
\newtheorem{prop}[thm]{Proposition}
\newtheorem{claim}{Claim}

\theoremstyle{definition}

\theoremstyle{remark}
\newtheorem{rem}{Remark}[section]

\begin{document}

\title{Stability of ideal lattices from quadratic number fields}
\author{Lenny Fukshansky}\thanks{The author was partially supported by the NSA Young Investigator Grant \#1210223 and a collaboration grant from the Simons Foundation (\#208969 to Lenny Fukshansky).}

\address{Department of Mathematics, 850 Columbia Avenue, Claremont McKenna College, Claremont, CA 91711}
\email{lenny@cmc.edu}


\begin{abstract}
We study semi-stable ideal lattices coming from real quadratic number fields. Specifically, we demonstrate infinite families of semi-stable and unstable ideal lattices of trace type, establishing explicit conditions on the canonical basis of an ideal that ensure stability; in particular, our result implies that an ideal lattice of trace type coming from a real quadratic field is semi-stable with positive probability. We also briefly discuss the connection between stability and well-roundedness of Euclidean lattices.
\end{abstract}

\maketitle

\def\A{{\mathcal A}}
\def\AA{{\mathfrak A}}
\def\B{{\mathcal B}}
\def\C{{\mathcal C}}
\def\D{{\mathcal D}}
\def\EE{{\mathfrak E}}
\def\F{{\mathcal F}}
\def\x{{\mathcal H}}
\def\I{{\mathcal I}}
\def\II{{\mathfrak I}}
\def\J{{\mathcal J}}
\def\K{{\mathcal K}}
\def\kk{{\mathfrak K}}
\def\L{{\mathcal L}}
\def\LL{{\mathfrak L}}
\def\M{{\mathcal M}}
\def\mm{{\mathfrak m}}
\def\MM{{\mathfrak M}}
\def\N{{\mathcal N}}
\def\O{{\mathcal O}}
\def\OO{{\mathfrak O}}
\def\PP{{\mathfrak P}}
\def\R{{\mathcal R}}
\def\PNR{{\mathcal P_N(\real)}}
\def\PMNR{{\mathcal P^M_N(\real)}}
\def\PdNR{{\mathcal P^d_N(\real)}}
\def\s{{\mathcal S}}
\def\V{{\mathcal V}}
\def\X{{\mathcal X}}
\def\Y{{\mathcal Y}}
\def\Z{{\mathcal Z}}
\def\H{{\mathcal H}}
\def\cee{{\mathbb C}}
\def\Nn{{\mathbb N}}
\def\pee{{\mathbb P}}
\def\que{{\mathbb Q}}
\def\QQ{{\mathbb Q}}
\def\real{{\mathbb R}}
\def\RR{{\mathbb R}}
\def\zed{{\mathbb Z}}
\def\ZZ{{\mathbb Z}}
\def\aaa{{\mathbb A}}
\def\ff{{\mathbb F}}
\def\HDelta{{\it \Delta}}
\def\kk{{\mathfrak K}}
\def\qbar{{\overline{\mathbb Q}}}
\def\kbar{{\overline{K}}}
\def\ybar{{\overline{Y}}}
\def\kkbar{{\overline{\mathfrak K}}}
\def\ubar{{\overline{U}}}
\def\eps{{\varepsilon}}
\def\ahat{{\hat \alpha}}
\def\bhat{{\hat \beta}}
\def\gt{{\tilde \gamma}}
\def\h{{\tfrac12}}
\def\be{{\boldsymbol e}}
\def\bei{{\boldsymbol e_i}}
\def\bc{{\boldsymbol c}}
\def\bm{{\boldsymbol m}}
\def\bk{{\boldsymbol k}}
\def\bi{{\boldsymbol i}}
\def\bl{{\boldsymbol l}}
\def\bq{{\boldsymbol q}}
\def\bu{{\boldsymbol u}}
\def\bt{{\boldsymbol t}}
\def\bs{{\boldsymbol s}}
\def\bv{{\boldsymbol v}}
\def\bw{{\boldsymbol w}}
\def\bx{{\boldsymbol x}}
\def\bX{{\boldsymbol X}}
\def\bz{{\boldsymbol z}}
\def\bwy{{\boldsymbol y}}
\def\bY{{\boldsymbol Y}}
\def\bL{{\boldsymbol L}}
\def\ba{{\boldsymbol a}}
\def\bb{{\boldsymbol b}}
\def\bet{{\boldsymbol\eta}}
\def\bxi{{\boldsymbol\xi}}
\def\bo{{\boldsymbol 0}}
\def\bone{{\boldsymbol 1}}
\def\bol{{\boldsymbol 1}_L}
\def\ep{\varepsilon}
\def\p{\boldsymbol\varphi}
\def\q{\boldsymbol\psi}
\def\rank{\operatorname{rank}}
\def\aut{\operatorname{Aut}}
\def\lcm{\operatorname{lcm}}
\def\sgn{\operatorname{sgn}}
\def\spn{\operatorname{span}}
\def\md{\operatorname{mod}}
\def\Norm{\operatorname{Norm}}
\def\dim{\operatorname{dim}}
\def\det{\operatorname{det}}
\def\Vol{\operatorname{Vol}}
\def\rk{\operatorname{rk}}
\def\ord{\operatorname{ord}}
\def\ker{\operatorname{ker}}
\def\div{\operatorname{div}}
\def\Gal{\operatorname{Gal}}
\def\GL{\operatorname{GL}}
\def\SNR{\operatorname{SNR}}
\def\WR{\operatorname{WR}}
\def\IWR{\operatorname{IWR}}
\def\scg{\operatorname{\left< \Gamma \right>}}
\def\swrh{\operatorname{Sim_{WR}(\Lambda_h)}}
\def\ch{\operatorname{C_h}}
\def\cht{\operatorname{C_h(\theta)}}
\def\scgt{\operatorname{\left< \Gamma_{\theta} \right>}}
\def\scgmn{\operatorname{\left< \Gamma_{m,n} \right>}}
\def\gat{\operatorname{\Omega_{\theta}}}
\def\mn{\operatorname{mn}}
\def\disc{\operatorname{disc}}
\def\tr{\operatorname{\tau_{\rho'}}}

\section{Introduction and statement of results}
\label{intro}

Let $\Lambda \subset \real^n$ be a lattice of rank $n \geq 2$. For each $1 \leq i \leq n$, the $i$-th successive minimum of $\Lambda$ is defined as
$$\lambda_i = \min \left\{ \lambda \in \real_{>0} : \dim \left( \spn_{\real} \left\{ \Lambda \cap B_n(\lambda) \right\} \right) \geq i \right\},$$
where $B_n(\lambda)$ is a closed ball of radius $\lambda$ centered at the origin in~$\real^n$. Then clearly
\begin{equation}
\label{s-min}
\lambda_1 \leq \dots \leq \lambda_n,
\end{equation}
and we say that $\Lambda$ is well-rounded (abbreviated WR) if there is equality throughout in~\eqref{s-min}. Two lattices $\Lambda$ and $\Omega$ are said to be similar, written $\Lambda \sim \Omega$, if there exists a positive real number $\gamma$ and an $n \times n$ real orthogonal matrix~$U$ such that~$\Lambda = \gamma U \Omega$. It is easy to see that ratios of successive minima, and hence well-roundedness, are preserved under similarity.

On the other hand, the lattice $\Lambda$ is called semi-stable if for each sublattice $\Omega \subseteq \Lambda$,
\begin{equation}
\label{semi_def}
\det(\Lambda)^{1/\rk(\Lambda)} \leq \det(\Omega)^{1/\rk(\Omega)}.
\end{equation}
For instance, when $\rk(\Lambda)=2$ the defining inequality~\eqref{semi_def} can be restated as
\begin{equation}
\label{semi_def-2}
\lambda_1 \geq \det(\Lambda)^{1/2},
\end{equation}
since for each sublattice $\Omega = \spn_{\zed} \left\{ \bz \right\} \subset \Lambda$ of rank~1, $\det(\Omega) = \|\bz\| \geq \lambda_1$. Semi-stability, the same as well-roundedness, is preserved under similarity. If a lattice is not semi-stable, we will say that it is unstable.

The notion of semi-stability was originally introduced by Stuhler~\cite{stuhler} in the context of reduction theory and later used by Grayson~\cite{grayson} in the study of arithmetic subgroups of semi-simple algebraic groups (see also~\cite{casselman} for an excellent survey of Stuhler's and Grayson's work). As indicated in~\cite{andre}, semi-stability heuristically means that the successive minima are not far from each other (see~\cite{borek} for a detailed investigation of this connection), i.e., inequality~\eqref{s-min} is not far from equality.  As a first observation, however, we note that the converse is not true; in other words, successive minima being close to each other does not necessarily imply stability. Specifically, we prove the following lemma.

\begin{lem} \label{wr_stable} All WR full-rank lattices in~$\real^2$ are semi-stable. On the other hand, for each $n \geq 3$ there exist infinitely many similarity classes of unstable WR lattices of rank~$n$ in $\real^n$.
\end{lem}

\proof
First suppose that $\Lambda \subset \real^2$ is WR. Then there exists a basis $\bx_1,\bx_2$ for $\Lambda$ consisting of vectors corresponding to successive minima, i.e.
$$\lambda_1 = \|\bx_1\| = \|\bx_2\| = \lambda_2.$$ 
Let $\theta$ be the angle between these vectors, then
$$\det(\Lambda) = \|\bx_1\| \|\bx_2\| \sin \theta = \lambda_1^2 \sin \theta \leq \lambda_1^2.$$
and so $\Lambda$ is semi-stable by~\eqref{semi_def-2}. This shows that all WR lattices in~$\real^2$ are semi-stable.

Next suppose $n \geq3$ and let $\be_1,\dots,\be_n$ be the standard basis vectors in~$\real^n$. We construct a family of examples of WR lattices of rank $n$ in $\real^n$, which are unstable. From our simple construction, it becomes immediately clear that many other such examples are possible. Let $\theta \in [\pi/3,\pi/2)$, and let 
$$\bx_{\theta} = \cos \theta \be_1 + \sin \theta \be_2,$$
and define
$$\Lambda_{\theta} = \spn_{\zed} \left\{ \be_1,\bx_{\theta},\be_3,\dots,\be_n \right\}.$$
It is easy to see that $\Lambda_{\theta}$ is WR with
$$\lambda_1= \dots = \lambda_n = 1,$$
where $\be_1,\bx_{\theta},\be_3,\dots,\be_n$ are the vectors corresponding to successive minima. Consider a sublattice $\Omega_{\theta} = \spn_{\zed} \left\{ \be_1,\bx_{\theta} \right\} \subset \Lambda$ of rank 2, and notice that
$$\det(\Lambda_{\theta})^{1/n} = \left( \sin \theta \right)^{1/n} > \left( \sin \theta \right)^{1/2} = \det(\Omega_{\theta})^{1/2},$$
since $\sqrt{3}/2 \leq \sin \theta < 1$. Hence $\Lambda_{\theta}$ is unstable, and two such lattices $\Lambda_{\theta_1}$ and $\Lambda_{\theta_2}$ are similar if and only if $\theta_1 = \theta_2$.
\endproof

\begin{rem} \label{kim} A particularly important subclass of WR lattices are perfect lattices, which figure prominently as potential candidates for extremum points of the sphere packing density function on the space of lattices, as well as in other related optimization problems. Y. Kim recently showed~\cite{ykim} that, while all perfect lattices in dimensions $\leq 7$ are semi-stable, there exists one 8-dimensional perfect lattice which is not semi-stable.
\end{rem}

In~\cite{andre}, the author remarks that, while semi-stable lattices have been investigated in several arithmetic and geometric contexts, they have not yet been seriously studied in the scope of classical lattice theory. A goal of this note is to partially remedy this situation. One important construction widely used in lattice theory is that of ideal lattices coming from number fields. Ideal lattices have been extensively studied in a series of papers by Eva Bayer-Fluckiger and her co-authors in the 1990's and 2000's (see, for instance, \cite{bayer1}, \cite{bayer2}, \cite{bayer_nebe}). Here we consider a restricted notion of ideal lattices coming from quadratic number fields, called ideal lattices of trace type. Let $K$ be a quadratic number field, and let us write $\O_K$ for its ring of integers. Then $K=\que(\sqrt{D})$ (real quadratic) or $K=\que(\sqrt{-D})$ (imaginary quadratic), where $D$ is a positive squarefree integer. The embeddings $\sigma_1, \sigma_2 : K \to \cee$ can be used  to define the standard Minkowski embedding $\sigma_K$ of $K$ into $\real^2$: if $K=\que(\sqrt{D})$, then $\sigma_K : K \to \real^2$ is given by $\sigma_K = (\sigma_1,\sigma_2)$; if $K=\que(\sqrt{-D})$, then $\sigma_2=\overline{\sigma_1}$, and $\sigma_K=(\Re(\sigma_1), \Im(\sigma_1))$, where $\Re$ and $\Im$ stand for real and imaginary parts, respectively. Each nonzero ideal $I \subseteq \O_K$ becomes a lattice of full rank in $\real^2$ under this embedding, which we will denote by $\Lambda_K(I) := \sigma_K(I)$. These are the ideal lattices we consider.

WR ideal lattices were studied in~\cite{lf:petersen} and~\cite{wr_ideal-2}, where in particular  it was shown that a positive proportion of quadratic number fields contain ideals giving rise to WR lattices. In view of Lemma~\ref{wr_stable}, it is interesting to understand which ideal lattices coming from quadratic number fields are semi-stable. An inequality connecting successive minima of an ideal lattice and the norm of its corresponding ideal~$I$ in the ring of integers of a fixed number field~$K$ follows from Lemma~3.2 of~\cite{lf:petersen}:
\begin{equation}
\label{s_min-to-norm}
\lambda_1(\Lambda_K(I))^2 \geq (r_1+r_2) \Nn(I)^{\frac{1}{r_1+r_2}}.
\end{equation}
Here $r_1$ is the number of real embeddings and~$r_2$ is the number of pairs of complex conjugate embeddings of~$K$; $\Nn(I)$ stands for the norm of the ideal~$I$ in~$\O_K$. A direct adaptation of Lemma~2 on p.115 of~\cite{lang} implies that 
\begin{equation}
\label{det_norm}
\det(\Lambda_K(I)) = 2^{-r_2} |\Delta_K|^{\frac{1}{2}} \Nn(I),
\end{equation}
where $\Delta_K$ is the discriminant of $K$. 

In this note, we discuss the case of real quadratic fields. When $K$ is a real quadratic number field, $r_1=2$ and $r_2=0$, and so combining~\eqref{s_min-to-norm} with~\eqref{det_norm}, we only obtain
$$\lambda_1(\Lambda_K(I)) \geq \frac{\sqrt{2}}{|\Delta_K|^{1/8}}\ \det(\Lambda_K(I))^{1/4}.$$
Hence the situation is more complicated and requires more detailed analysis and additional notation. Let~$D > 1$ be a squarefree integer and let $K=\que(\sqrt{D})$. We have $\O_K=\zed[\delta]$, where
\begin{equation}
\label{delta}
\delta = \left\{ \begin{array}{ll}
- \sqrt{D} & \mbox{if $K=\que(\sqrt{D})$, $D \not\equiv 1 (\md 4)$} \\
\frac{1-\sqrt{D}}{2} & \mbox{if $K=\que(\sqrt{D})$, $D \equiv 1 (\md 4)$.}
\end{array}
\right.
\end{equation}
Now $I \subseteq \O_K$ is an ideal if and only if 
\begin{equation}
\label{I_abg}
I = I(a,b,g) :=  \{ ax + (b+g\delta)y : x,y \in \zed \},
\end{equation}
for some $a,b,g \in \zed_{\geq 0}$ such that
\begin{equation}
\label{abg}
b < a,\ g \mid a,b,\text{ and } ag \mid \Nn(b+g\delta).
\end{equation}
Such integral basis $a,b+g\delta$ is unique for each ideal $I$ and is called the canonical basis for~$I$ (see Section~6.3 of~\cite{buell} for details). In Section~\ref{real} we prove the following result.

\begin{thm} \label{real_quad} Let $K = \que(\sqrt{D})$ be a real quadratic number field. Then there exist infinitely many ideals $I \subseteq \O_K$ for which the corresponding ideal lattice~$\Lambda_K(I)$ is semi-stable, as well as infinitely many such ideals with the corresponding lattice unstable. Specifically, let $\gamma \in \real_{>0}$ and define the functions
$$u_{\gamma}(b) = \left\{ \begin{array}{ll}
\frac{\gamma(2b+1)}{2} & \mbox{if $D \equiv 1 (\md 4)$} \\
\gamma b & \mbox{if $D \not\equiv 1 (\md 4)$,}
\end{array}
\right.$$
$$v(b) = \left\{ \begin{array}{ll}
\frac{(2b+1)^2+D}{2\sqrt{D}} & \mbox{if $D \equiv 1 (\md 4)$} \\
\frac{b^2+D}{\sqrt{D}} & \mbox{if $D \not\equiv 1 (\md 4)$,}
\end{array}
\right.$$
$$h(b) = \left\{ \begin{array}{ll}
\frac{(2b+1)^2-D}{2} & \mbox{if $D \equiv 1 (\md 4)$} \\
b^2-D & \mbox{if $D \not\equiv 1 (\md 4)$.}
\end{array}
\right.$$
Then there exists an absolute constant $\gamma > 1$ such that if
\begin{equation}
\label{stable-1}
u_{\gamma}(b) \leq a \leq v(b),
\end{equation}
then the lattice $\Lambda_K(I(a,b,g))$ is semi-stable for every triple $a,b,g$ satisfying~\eqref{abg}. On the other hand, if
\begin{equation}
\label{nonstable-1}
v(b) < a \leq h(b),
\end{equation}
then the lattice $\Lambda_K(I(a,b,g))$ is unstable for every triple $a,b,g$ satisfying~\eqref{abg}. 
\end{thm}

\noindent
In fact, Remark~\ref{proportion} below shows that the probability of an arbitrary ideal lattice~$\Lambda_{\que(\sqrt{D})}(I(a,b,g))$ being semi-stable is positive (specifically, the probability is at least $1/\gamma$ as $b \to \infty$).
\smallskip

In Section~\ref{auxil} we prove a technical lemma on distribution of divisors of integers of the form $x^2-D$, which is useful to us later in our main argument. In Section~\ref{real_lemmas} we establish Proposition~\ref{ab_stable}, which is the core of our argument. Finally, we use this proposition in Section~\ref{real} to prove Theorem~\ref{real_quad}. We are now ready to proceed.
\bigskip

\section{A divisor lemma}
\label{auxil}

In this section we make an observation on the finiteness of the set of integers of the form $x^2 \pm D$ which have divisors in small intervals around their square root. This result is later used in the proof of Theorem~\ref{real_quad}. The proof of this lemma was suggested to me by Florian Luca.

\begin{lem} \label{fin_div} Let $|D| > 1$ be a squarefree integer and $0 < \eps < 1/2$ a real number. Then the set
$$\left\{ x \in \zed_{>0} : \exists\ b \mid x^2-D \text{ such that } x < b \leq x+x^{1/2-\eps} \right\}$$
is finite.
\end{lem}

\proof
Since there are only finitely many positive integers less than any fixed constant, we can assume without loss of generality that
$$x > \max \left\{ |D|, 2^{1/\eps} \right\}.$$
Let us write $x^2-D = bd$, where $b \in (x, x+x^{1/2-\eps}]$, then $d \in [x-x^{1/2-\eps},x)$. Notice that $b=d+a$, where $a \in [0,2x^{1/2-\eps}]$. Therefore
$$x^2-D = d(d+a) = d^2 + 2d \frac{a}{2} + \left( \frac{a}{2} \right)^2 -  \left( \frac{a}{2} \right)^2 = \left( d+\frac{a}{2} \right)^2 -  \left( \frac{a}{2} \right)^2,$$
and therefore
$$(2x)^2 = (2d+a)^2 + (4D- a^2),$$
meaning that
$$(2x-(2d+a))(2x+(2d+a)) = 4D - a^2.$$
Taking absolute values, we see that the left hand side cannot be equal to zero; since $|D| > 1$, the assumption that $4D - a^2 = 0$ would imply that $D = (a/2)^2 > 1$, which would contradict $D$ being squarefree. Since $2x-(2d+a)$ is an integer, $|2x-(2d+a)| \geq 1$, which means that
$$|(2x-(2d+a))(2x+(2d+a))| \geq 2x+(2d+a) > 2x.$$
On the other hand,
$$|4D - a^2| \leq 4|D|+a^2 < 4|D| + 4x^{1-2\eps},$$
and so we have
$$2x < 4x^{1-2\eps} + 4|D|.$$
Therefore, since $x > 2^{1/\eps}$,
$$x < 2x(1 - 2x^{-2\eps}) < 4|D|,$$
meaning that there are at most $4|D|$ such integers $x$.
\endproof
\bigskip

\section{Lemmas on stability of some planar lattices}
\label{real_lemmas}

Our goal here is to develop a collection of lemmas that will allow us to treat ideal lattices coming from any real quadratic number field simultaneously. Throughout this section, let~$D > 1$ be fixed a squarefree integer. For each pair of integers $(a,b)$ such that
\begin{equation}
\label{ab1}
0 < b < a,\ a \mid b^2-D,
\end{equation}
define the lattice
\begin{equation}
\label{Lab}
\Lambda(a,b) = \begin{pmatrix} a & b-\sqrt{D} \\ a & b+\sqrt{D} \end{pmatrix} \zed^2.
\end{equation}
We want to understand for which pairs $(a,b)$ satisfying~\eqref{ab1} the corresponding lattice $\Lambda(a,b)$ is semi-stable. Let
$$S(D) = \left\{ (a,b) \in \zed^2 : (a,b) \text{ satisfies \eqref{ab1}} \right\}.$$
We prove the following result.

\begin{prop} \label{ab_stable} For infinitely many pairs $(a,b) \in S(D)$, the corresponding lattice $\Lambda(a,b)$ is semi-stable, and for infinitely many pairs it is unstable. Specifically, there exists an absolute constant $\gamma > 1$ such that if
\begin{equation}
\label{fm7}
\gamma b \leq a \leq \frac{b^2+D}{\sqrt{D}},
\end{equation}
then the lattice $\Lambda(a,b)$ is semi-stable. On the other hand, if
\begin{equation}
\label{fm8}
\frac{b^2+D}{\sqrt{D}} < a \leq b^2-D,
\end{equation}
then the lattice $\Lambda(a,b)$ is unstable.
\end{prop}

To establish Proposition~\ref{ab_stable}, notice that for each $(a,b) \in S(D)$, $\det(\Lambda(a,b)) = 2a\sqrt{D}$, and so $\Lambda(a,b)$ is semi-stable if and only if
$$\lambda_1(\Lambda(a,b))^2 \geq 2a \sqrt{D}.$$
The norm form of $\Lambda(a,b)$ corresponding to the choice of basis as in~\eqref{Lab} is
$$Q(x,y) = Q_{(a,b)}(x,y) := 2(x a + y b)^2 + 2y^2 D,$$
then
$$\lambda_1^2 = \min \left\{ Q(x,y) : (x,y) \in \zed^2 \setminus \{ (0,0) \} \right\}.$$
Let $(\alpha,\beta) \in \zed^2$ be a point at which this minimum is achieved, i.e.,
$$Q(\alpha,\beta) = 2 \min \left\{ (x a + y b)^2 + y^2 D : (x,y) \in \zed^2 \setminus \{ (0,0) \} \right\},$$
then $\gcd(\alpha,\beta)=1$, and semi-stability is equivalent to the inequality
\begin{equation}
\label{form_min}
Q(\alpha,\beta) \geq 2 a \sqrt{D}.
\end{equation}

\begin{lem} \label{123} $(\alpha,\beta)$, the minimum of $Q(x,y)$ falls into one of the following three categories:
\begin{enumerate}
\item $(\alpha,\beta) = (1,0)$,
\item $(\alpha,\beta) = (0,1)$,
\item $0 < \alpha \leq b,\ \alpha \leq |\beta| \leq a,\ \beta < 0$.
\end{enumerate}
\end{lem}

\proof
Assume (I) and (II) do not hold, which means that $\alpha \beta \neq 0$. Then $\alpha \beta < 0$, since otherwise 
$$Q(\alpha,\beta) \geq 2 (a+b)^2 + 2D > 2(b^2+D) = Q(0,1).$$
Hence we can assume without loss of generality that $\beta < 0$, since $Q(\alpha,\beta) = Q(-\alpha,-\beta)$. If $|\beta| > a$, then 
$$Q(\alpha,\beta) > 2a^2D > 2a^2 = Q(1,0).$$
Now consider 
$$f(\alpha) = Q(\alpha,\beta) = 2(\alpha a - |\beta| b)^2 + 2\beta^2 D$$
as a function of $\alpha$. Notice that it is increasing when $\alpha > |\beta| b/a$. Since $|\beta| \leq a$, $\alpha > |\beta| b/a$ when $\alpha > b$, meaning that $Q(\alpha,\beta)$ cannot achieve its minimum for such values of~$\alpha$. Finally, assume that $\alpha > |\beta|$ and recall that $a > b$. Then
$$Q(\alpha,\beta) = 2(\alpha a - |\beta| b)^2 + 2\beta^2 D \geq 2(2a-b)^2 + 2D = 8(a^2-ab) + 2(b^2+D) > Q(0,1).$$
Hence we established that the inequalities (III) hold, which proves the lemma.
\endproof

Let us define three sets of pairs $(a,b) \in S(D)$, corresponding to each of the three cases above:
$$S_1 = S_1(D) := \left\{ (a,b) \in S(D) : (\alpha,\beta) \text{ is as in (I)} \right\},$$
$$S_2 = S_2(D) := \left\{ (a,b) \in S(D) : (\alpha,\beta) \text{ is as in (II)} \right\},$$
$$S_3 = S_3(D) := \left\{ (a,b) \in S(D) : (\alpha,\beta) \text{ is as in (III)} \right\}.$$

\noindent
We can write $a = C(b^2-D)$ for some $C \in \real_{>0}$, $b/(b^2-D) < C \leq 1$. Then $\Lambda(a,b)$ is semi-stable if and only if $Q(\alpha,\beta) \geq 2 C(b^2-D) \sqrt{D}$, which is equivalent to
\begin{equation}
\label{fm1}
C \left( \alpha^2 C(b^2-D) + 2\alpha \beta b - \sqrt{D} \right) \geq - \frac{\beta^2(b^2+D)}{b^2-D}.
\end{equation}
The right hand side of~\eqref{fm1} is always non-positive and $C > 0$.
\smallskip

\begin{lem} \label{lem-S1} The set $S_1$ is finite, and the lattice $\Lambda(a,b)$ is semi-stable for every pair $(a,b) \in S_1$ with $b > \sqrt{D}$.
\end{lem}

\proof
Let $(a,b) \in S_1$ with $b > \sqrt{D}$, then $\beta = 0$, $\alpha = 1$ and~\eqref{fm1} holds for all values of~$C$. Hence the lattice $\Lambda(a,b)$ is semi-stable.

Now we show that $S_1$ is finite. Notice that for each $(a,b) \in S_1$,
$$\frac{1}{2} Q(1,0) = C^2 (b^2-D)^2 \leq \frac{1}{2} Q(0,1) = b^2+D,$$
and so $C \leq \frac{\sqrt{b^2+D}}{b^2-D}$, which means that
$$b < a \leq \sqrt{b^2+D} < b + \sqrt{D}.$$
Therefore $b$ is an integer such that $b^2-D$ has a divisor $a \in (b, b + \sqrt{D})$, and clearly $\sqrt{D} < b^{1/2-\eps}$ for any $\eps > 0$ for all but finitely many $b$. There are only finitely many such integers $b$ by Lemma~\ref{fin_div}, and so the set of such pairs $(a,b)$ is finite, since $a$ is bounded by~$\sqrt{b^2+D}$.
\endproof

\begin{lem} \label{lem-S2} Let $(a,b) \in S_2$ and $a = C(b^2-D)$ as above. Then $\Lambda(a,b)$ is semi-stable if and only if $C \leq \frac{b^2+D}{(b^2-D)\sqrt{D}}$.
\end{lem}

\proof
Suppose $\alpha = 0$, then $\beta = 1$, then~\eqref{fm1} holds if and only if
\begin{equation}
\label{fm2}
C \leq \frac{b^2+D}{(b^2-D)\sqrt{D}}.
\end{equation}
\endproof

\begin{lem} \label{lem-S3} Let $(a,b) \in S_3$ and $a = C(b^2-D)$ as above. There exists an absolute real constant $\gamma > 1$ such that if $C \geq \frac{\gamma b}{b^2-D}$, then $\Lambda(a,b)$ is semi-stable.
\end{lem}

\proof
If the set $S_3$ is finite, there is nothing to prove, so assume it is infinite. Let 
$$S_3' = \left\{ b \in \zed_{>0} : \exists\ a \in \zed_{>0} \text{ such that } (a,b) \in S_3 \right\}.$$
In the asymptotic argument below, when we consider $b$ getting large or tending to infinity, we always mean that $b$ stays in $S_3'$ and $a = C(b^2-D)$ is such that $(a,b) \in S_3$.

For each $(a,b) \in S_3$, the corresponding $\alpha,\beta \neq 0$ are such that $\beta < 0 < \alpha \leq |\beta|$.  The inequality~\eqref{fm1} certainly holds when
\begin{equation}
\label{fm3}
\alpha^2 C(b^2-D) + 2\alpha \beta b - \sqrt{D}  \geq 0,
\end{equation}
which is true whenever
\begin{equation}
\label{fm4}
C \geq \frac{2 \alpha |\beta| b + \sqrt{D}}{\alpha^2 (b^2-D)} = \left( \frac{|\beta|}{\alpha} \right) \left( \frac{2 b}{b^2-D} \right) + \frac{\sqrt{D}}{\alpha^2 (b^2-D)}.
\end{equation}

\begin{claim} \label{clm1} There exists an absolute constant $\rho$ so that $1 \leq \frac{|\beta|}{\alpha} \leq \rho$ for all~$(a,b) \in S_3$.
\end{claim}

\proof
Suppose not, then there exists some monotone increasing unbounded real-valued function $f(b)$ such that
\begin{equation}
\label{rho}
\liminf_{b \to \infty} \frac{|\beta|}{\alpha f(b)} = 1.
\end{equation}
Hence we can assume that there exists an infinite subsequence of positive integers $b$ for which $\beta \sim - \alpha f(b)$ as $b \to \infty$. Then for all sufficiently large $b$,
\begin{eqnarray}
\label{fm5}
& & \frac{1}{2} \min \left\{ Q(x,y) : (x,y) \in \zed^2 \setminus {(0,0)} \right\} \nonumber \\
& = & \frac{1}{2} Q(\alpha,\beta) \sim \left( \alpha C (b^2-D) - \alpha b f(b) \right)^2 + \alpha^2 f(b)^2 D \nonumber \\
& = &  \alpha^2 f(b)^2 \left(b^2 \left( \frac{C(b^2-D)}{b f(b)} - 1 \right)^2 + D \right) > b^2+D = \frac{1}{2} Q(0,1),
\end{eqnarray}
unless $\frac{C(b^2-D)}{b f(b)} \to 1$ as $b \to \infty$. Suppose this is the case, then
$$\frac{C(b^2-D)}{b f(b)} = \left( \frac{|\beta|}{\alpha f(b)} \right) \frac{\alpha C(b^2-D)}{ |\beta| b} \to 1,$$
and since $|\beta|/\alpha f(b) \to 1$ and $a = C(b^2-D)$, we have $\frac{a}{b} \times \frac{\alpha}{|\beta|} \to 1$ as $b \to \infty$. Since $a,b,\alpha,\beta$ are integers, we must have
\begin{equation}
\label{aabb}
\frac{a}{b} = \frac{|\beta|}{\alpha}
\end{equation}
for all sufficiently large $b$. Since $\alpha$ and $\beta$ are relatively prime, we must have $\alpha=b/d$, $\beta=-a/d$, where $d=\gcd(a,b) \mid D$ by~\eqref{ab1}. Then
$$\frac{1}{2} Q(\alpha,\beta) = \frac{a^2 D}{d^2} \leq b^2+D = \frac{1}{2} Q(0,1),$$
and so
\begin{equation}
\label{aabb-1}
\frac{a}{b} \leq d \sqrt{ \frac{1}{D} + \frac{1}{b^2} } < d \sqrt{2} \leq D\sqrt{2}.
\end{equation}
Now \eqref{aabb-1} combined with \eqref{aabb} implies that $|\beta|/\alpha \leq D\sqrt{2}$. This completes the proof.
\endproof

Thus we conclude that $|\beta|/\alpha \leq \rho$ for all~$b \in S_3'$. Then~\eqref{fm4} implies that for all ~$b \in S_3'$, if
\begin{equation}
\label{fm6}
C \geq \frac{2 \rho b}{b^2-D} + \frac{\sqrt{D}}{\alpha^2 (b^2-D)},
\end{equation}
then the lattice $\Lambda(a,b)$ is semi-stable. In other words, there exists some real constant $\gamma > \rho \geq 1$ such that whenever $a=C(b^2-D)$ for $C \in [\gamma b /(b^2-D), 1]$ so that $(a,b) \in S_3$, the lattice $\Lambda(a,b)$ is semi-stable.
\endproof
\smallskip

\proof[Proof of Proposition~\ref{ab_stable}]
Let $\gamma$ be the constant as in the statement of Lemma~\ref{lem-S3}. First, let $(a,b) \in S(D)$ as above with $b > \sqrt{D}$, and assume that~\eqref{fm7} is satisfied. Notice that $(a,b)$ is either in $S_1$, $S_2$, or $S_3$. Then the result follows by combining Lemmas~\ref{lem-S1}, \ref{lem-S2}, and~\ref{lem-S3}.

Next, assume that~\eqref{fm8} holds, then
$$\det(\Lambda(a,b)) = 2a \sqrt{D} > 2(b^2+D) = Q(0,1) \geq \lambda_1(\Lambda(a,b))^2,$$
and so $\Lambda(a,b)$ is unstable.

To construct an infinite family of pairs $(a,b) \in S(D)$ giving rise to unstable lattices, simply take $a=b^2-D$ for each integer $b > \sqrt{\frac{D+D\sqrt{D}}{\sqrt{D}-1}}$; the resulting lattice is unstable since~\eqref{fm8} is satisfied.

On the other hand, for each $m \in \zed_{>0}$ let $b=mD$ and take $a = \frac{b^2-D}{D} = m^2D-1$. Let $\gamma$ be the constant as in the statement of Lemma~\ref{lem-S3}. For each $m \geq \frac{\gamma D + \sqrt{\gamma^2 D^2 + 4D}}{2D}$, the inequality~\eqref{fm7} is satisfied,  and hence the resulting lattice is semi-stable by the argument above.
\endproof
\smallskip

\begin{rem} \label{small_a} In the argument above, we constructed a family of unstable lattices $\Lambda(a,b)$ with $a$ large comparing to $b$. On the other hand, there also exist unstable lattices $\Lambda(a,b)$ with $a$ close to $b$. For instance, let $D=13$ and consider the pair $(a,b)=(276,259) \in S(D)$. Then
$$\lambda_1^2 \leq Q_{(a,b)}(1,-1) = 604 < 2a\sqrt{D} = 552\sqrt{13} = \det(\Lambda(a,b)),$$
and so the lattice $\Lambda(276,259)$ is unstable.
\end{rem}

\bigskip

\section{The case of real quadratic number fields}
\label{real}

In this section we prove Theorem~\ref{real_quad}. Let $D>1$ be a squarefree integer, $K=\que(\sqrt{D})$, integers $a,b,g \geq 0$ satisfying~\eqref{abg}, and the ideal $I=I(a,b,g) \subseteq \O_K$ as in~\eqref{I_abg}. Then
\begin{equation}
\label{B1}
\Lambda_K(I) = \begin{pmatrix} a & b-g\sqrt{D} \\ a & b+g\sqrt{D} \end{pmatrix} \zed^2,
\end{equation}
if $D \not\equiv 1 (\md 4)$, and
\begin{equation}
\label{B2}
\Lambda_K(I) = \begin{pmatrix} a & \frac{2b+g}{2} - \frac{g\sqrt{D}}{2} \\ a &\frac{2b+g}{2} + \frac{g\sqrt{D}}{2} \end{pmatrix} \zed^2,
\end{equation}
if $D \equiv 1 (\md 4)$. Notice that $I=gI'$, where $I'$ has canonical basis $\frac{a}{g}, \frac{b}{g}+\delta$ and $\Lambda_K(I) \sim \Lambda_K(I')$. Hence we can assume without loss of generality that $g=1$.
\smallskip

First assume that $D \not\equiv 1 (\md 4)$, then
$$I = \{ ax + (b-\sqrt{D})y : x,y \in \zed \} \subseteq \O_K.$$
Here the pair $(a,b)$ satisfies the conditions of~\eqref{ab1} and $\Lambda_K(I) = \Lambda(a,b)$. The statement of Theorem~\ref{real_quad} in this case readily follows from Proposition~\ref{ab_stable}.
\smallskip

Now assume that $D \equiv 1 (\md 4)$, then
$$I = \left\{ ax + \left( \frac{2b + 1-\sqrt{D}}{2} \right) y  :  x,y \in \zed \right\} \subseteq \O_K,$$
where
\begin{equation}
\label{abg-2}
b < a,\ a \mid \frac{1}{4} \left( (2b+1)^2 - D \right),
\end{equation}
and
\begin{equation}
\label{basis-2}
\Lambda_K(I) = \begin{pmatrix} a & \frac{2b+1}{2} - \frac{\sqrt{D}}{2} \\ a &\frac{2b+1}{2} + \frac{\sqrt{D}}{2} \end{pmatrix} \zed^2.
\end{equation}
Let $a_1=2a$, $b_1=2b+1$, and notice that the pair $(a_1,b_1)$ satisfies the conditions of~\eqref{ab1} and $\Lambda_K(I) = \frac{1}{2} \Lambda(a_1,b_1)$. Observe that $\Lambda_K(I)$ is semi-stable if and only if $\Lambda(a_1,b_1)$ is semi-stable, and hence the statement of Theorem~\ref{real_quad} in this case again follows from Proposition~\ref{ab_stable}.
\bigskip

\begin{rem} \label{proportion} In fact, Theorem~\ref{real_quad} implies that an arbitrary ideal lattice~$\Lambda_{\que(\sqrt{D})}(I)$ is semi-stable with positive probability.

Indeed, for $x > 0$ let
\begin{equation}
M(D,x) =  \left\{ q \in \zed_{>0} : q < x, D \text{ is a quadratic residue} \md q \right\},
\end{equation}
then $q < x$ is in $M(D,x)$ if and only if $D$ is a quadratic residue modulo every prime dividing~$q$. Professor Gang Yu pointed out to me that an argument essentially identical to the proof of the main result of~\cite{james} shows that there exists a positive real constant $C(D)$ such that
\begin{equation}
\label{MDx}
|M(D,x)| \sim C(D) \left( \frac{x}{\sqrt{\log x}} \right)
\end{equation}
as $x \to \infty$ (the set $M(D,x)$ can also be compared to the set $S(x)=M(-1,x)$ in the definition of Landau-Ramanujan constant~\cite{landau}, where the same classical asymptotic emerges). Then~\eqref{MDx} implies that for $0 < k_1 < k_2 < 1$,
$$\left| M(D,x) \cap [k_1x,k_2x] \right| \sim C(D) (k_2-k_1) \left( \frac{x}{\sqrt{\log x}} \right),$$
and so
$$\lim_{x \to \infty} \frac{\left| M(D,x) \cap [k_1x,k_2x] \right|}{\left| M(D,x) \right|} = k_2-k_1,$$
which means that elements of $M(D,x)$ are equidistributed in subintervals of $[1,x]$. In other words, as $x \to \infty$, every subinterval $[k_1x,k_2x]$ with $0 < k_1 < k_2 < 1$ will contain a $(k_2-k_1)$-proportion of integers $q$ such that $D$ is quadratic residue modulo~$q$. This implies that probability of such a modulus $q$ to be in the interval $[k_1x,k_2x]$ tends to $k_2-k_1$ as $x \to \infty$.

Now let $K = \que(\sqrt{D})$, let $I = I(a,b,1) \subseteq \O_K$ be an ideal, and let
$$a_1 = \left\{ \begin{array}{ll}
a & \mbox{if $D \not\equiv 1 (\md 4)$} \\
2a & \mbox{if $D \equiv 1 (\md 4)$,}
\end{array}
\right.$$
and
$$b_1 = \left\{ \begin{array}{ll}
b & \mbox{if $D \not\equiv 1 (\md 4)$} \\
2b+1 & \mbox{if $D \equiv 1 (\md 4)$.}
\end{array}
\right.$$
Then $a_1 \mid b_1^2 - D$ and $b_1 < a_1 \leq b_1^2 -D$. Let $d_1 = (b_1^2-D)/a_1$, so $d_1 \mid b_1^2 - D$ and $1 \leq d_1 < b_1$. Theorem~\ref{real_quad} implies that if
\begin{equation}
\label{d_bnd}
\sqrt{D} \left( \frac{b_1^2-D}{b_1^2+D} \right) \leq d_1 \leq \frac{1}{\gamma} b_1 - \frac{D}{\gamma b_1},
\end{equation}
then the lattice $\Lambda_K(I)$ is semi-stable. In other words,~\eqref{d_bnd} implies that for each $\eps > 0$ there exists $B \in \real_{>0}$ such that for all $b_1 > B$, if 
\begin{equation}
\label{d_bnd-1}
d_1 \in \left[ \frac{\sqrt{D}}{b_1} b_1, \frac{1}{\gamma} b_1 - \eps \right],
\end{equation}
then $\Lambda_K(I)$ is semi-stable. Since $d_1 \in M(D,b_1)$, our argument above suggests that the probability of~\eqref{d_bnd-1} holding tends to $\frac{1}{\gamma}$ as $b_1 \to \infty$.
\end{rem}

{\bf Acknowledgment.} I would like to thank Professor Florian Luca for suggesting the proof of Lemma~\ref{fin_div}, as indicated above. I also thank Professors Gang Yu and David Speyer, whose comments were instrumental to the formulation of Remark~\ref{proportion}. Finally, I thank the referee for many useful suggestions which improved the quality of the paper.


\bigskip

\bibliographystyle{plain}  
\bibliography{semi-stable}    

\begin{thebibliography}{10}

\bibitem{andre}
Y.~Andr\'e.
\newblock On nef and semistable hermitian lattices, and their behaviour under
  tensor product.
\newblock {\em Tohoku Math. J. (2)}, 63(4):629--649, 2011.

\bibitem{bayer1}
E.~Bayer-Fluckiger.
\newblock Lattices and number fields.
\newblock Contemp. Math. 241, pages 69--84, 1999.

\bibitem{bayer2}
E.~Bayer-Fluckiger.
\newblock Ideal lattices.
\newblock In {\em A panorama of number theory or the view from Baker's garden
  (Zurich, 1999)}, pages 168--184. Cambridge Univ. Press, Cambridge, 2002.

\bibitem{bayer_nebe}
E.~Bayer-Fluckiger and G.~Nebe.
\newblock On the {E}uclidean minimum of some real number fields.
\newblock {\em J. Th\'eor. Nombres Bordeaux}, 17(2):437Ð454, 2005.

\bibitem{borek}
T.~Borek.
\newblock Successive minima and slopes of {H}ermitian vector bundles over
  number fields.
\newblock {\em J. Number Theory}, 113(2):380--388, 2005.

\bibitem{buell}
D.~A. Buell.
\newblock {\em Binary Quadratic Forms}.
\newblock Springer-Verlag, 1989.

\bibitem{casselman}
B.~Casselman.
\newblock Stability of lattices and the partition of arithmetic quotients.
\newblock {\em Asian J. Math.}, 8(4):607--637, 2004.

\bibitem{wr_ideal-2}
L.~Fukshansky, G.~Henshaw, P.~Liao, M.~Prince, X.~Sun, and S.~Whitehead.
\newblock On well-rounded ideal lattices, {II}.
\newblock {\em Int. J. Number Theory}, 9(1):139--154, 2013.

\bibitem{lf:petersen}
L.~Fukshansky and K.~Petersen.
\newblock On ideal well-rounded lattices.
\newblock {\em Int. J. Number Theory}, 8(1):189--206, 2012.

\bibitem{grayson}
D.~R. Grayson.
\newblock Reduction theory using semistability.
\newblock {\em Comment. Math. Helv.}, 59(4):600--634, 1984.

\bibitem{james}
R.~D. James.
\newblock The distribution of integers represented by quadratic forms.
\newblock {\em Amer. J. Math.}, 60(3):737--744, 1938.

\bibitem{ykim}
Y.~Kim.
\newblock On semistability of root lattices and perfect lattices.
\newblock {\em preprint, Univ. Illinois}, 2009.

\bibitem{lang}
S.~Lang.
\newblock {\em Algebraic Number Theory}.
\newblock Springer-Verlag, 1994.

\bibitem{stuhler}
U.~Stuhler.
\newblock Eine {B}emerkung zur {R}eduktionstheorie quadratischer {F}ormen.
\newblock {\em Arch. Math. (Basel)}, 27(6):604--610, 1976.

\bibitem{landau}
E.~W. Weisstein.
\newblock Landau-{R}amanujan constant. {F}rom {M}ath{W}orld -- a {W}olfram web
  resource.
\newblock \url{http://mathworld.wolfram.com/Landau-RamanujanConstant.html}.

\end{thebibliography}
\end{document}